\documentclass[journal=jacsat,manuscript=article]{achemso}

\usepackage[version=3]{mhchem} 
\usepackage{xr}
\usepackage{cleveref}
\usepackage{mathptmx}
\usepackage{amsfonts}

\author{Guangqi Wu}
\affiliation{Department of Chemical Engineering, Massachusetts Institute of Technology, 77 Massachusetts Avenue, Cambrdige, MA 02139, USA. }
\alsoaffiliation{Department of Electrical Engineering and Computer Science, Massachusetts Institute of Technology, 77 Massachusetts Avenue, Cambrdige, MA 02139, USA.}

\altaffiliation{Current address: Department of Chemistry, University of Oxford, 12 Mansfield Road, Oxford, OX1 3TA, United Kingdom.}
\alsoaffiliation{Equal contribution}
\author{Runzhong Wang}
\affiliation{Department of Chemical Engineering, Massachusetts Institute of Technology, 77 Massachusetts Avenue, Cambrdige, MA 02139, USA. }
\alsoaffiliation{Department of Electrical Engineering and Computer Science, Massachusetts Institute of Technology, 77 Massachusetts Avenue, Cambrdige, MA 02139, USA.}
\alsoaffiliation{Equal contribution}
\author{Connor W. Coley}
\affiliation{Department of Chemical Engineering, Massachusetts Institute of Technology, 77 Massachusetts Avenue, Cambrdige, MA 02139, USA. }
\alsoaffiliation{Department of Electrical Engineering and Computer Science, Massachusetts Institute of Technology, 77 Massachusetts Avenue, Cambrdige, MA 02139, USA.}
\email{ccoley@mit.edu}

\title[An \textsf{achemso} demo]
  {Optimization of Robotic Liquid Handling as a Capacitated Vehicle Routing Problem}

\abbreviations{IR,NMR,UV}

\begin{document}


\begin{abstract}
We present an optimization strategy to reduce the execution time of liquid handling operations in the context of an automated chemical laboratory. 
By formulating the task as a capacitated vehicle routing problem (CVRP), we leverage heuristic solvers traditionally used in logistics and transportation planning to optimize task execution times. As exemplified using an 8-channel pipette with individually controllable tips, our approach demonstrates robust optimization performance across different labware formats (e.g., well-plates, vial holders), achieving up to a 37\% reduction in execution time for randomly generated tasks compared to the baseline sorting method. We further apply the method to a real-world high-throughput materials discovery campaign and observe that 3 minutes of optimization time led to a reduction of 61 minutes in execution time compared to the best-performing sorting-based strategy. Our results highlight the potential for substantial improvements in throughput and efficiency in automated laboratories without any hardware modifications. This optimization strategy offers a practical and scalable solution to accelerate combinatorial experimentation in areas such as drug combination screening, reaction condition optimization, materials development, and formulation engineering.
\end{abstract}

\section{Introduction}
Liquid handling systems play a central role in modern lab automation by relieving researchers from repetitive and time-intensive tasks. With the integration of computational methods for experimental design, automated platforms have shown extraordinary promise in scientific discovery, particularly in areas of life science, chemistry, materials, and drug discovery\cite{koscher2023autonomous,dai2024autonomous,shields2021bayesian,abolhasani2023rise,tom2024self,chen2024navigating,wang2025autonomous,mehr2020universal,stein2019progress}. 
While advances in algorithms continue to improve our sampling from the design space\cite{du2024machine,tran2024design,sanchez2018inverse}, our ability to translate these designs into actionable experiments remains constrained by practical considerations of execution time. 
The ability to efficiently access larger design spaces within a limited time frame is crucial for accelerated discovery.

Combinatorial screening has seen renewed attention within the realms of drug and materials discovery, with applications spanning drug combinations, polymers, formulations, and battery materials\cite{jin2023rational,sun2013high,bannigan2023machine,tamasi2022machine,reis2021machine,noh2024integrated,ting2016high}. 
In these workflows, liquid handling plays a critical role in transferring material from stock solutions to each (combinatorial) mixture to be evaluated. 
With a sufficiently fast downstream assay (e.g., an optical measurement, direct injection mass spectrometry), the most time-consuming step in a combinatorial screen is liquid handling. 
As the number of potential components increases (both in terms of the number of distinct stock solutions and the number of distinct components that might be included in each mixture), execution bottlenecks become more severe. In our own experience, combinatorial liquid handling involving approximately 350 transfers from one 96-well plate to another on a Tecan Evo 200 liquid handler requires upwards of an hour to execute.

\begin{figure*}
 \centering
 \includegraphics[height=8cm]{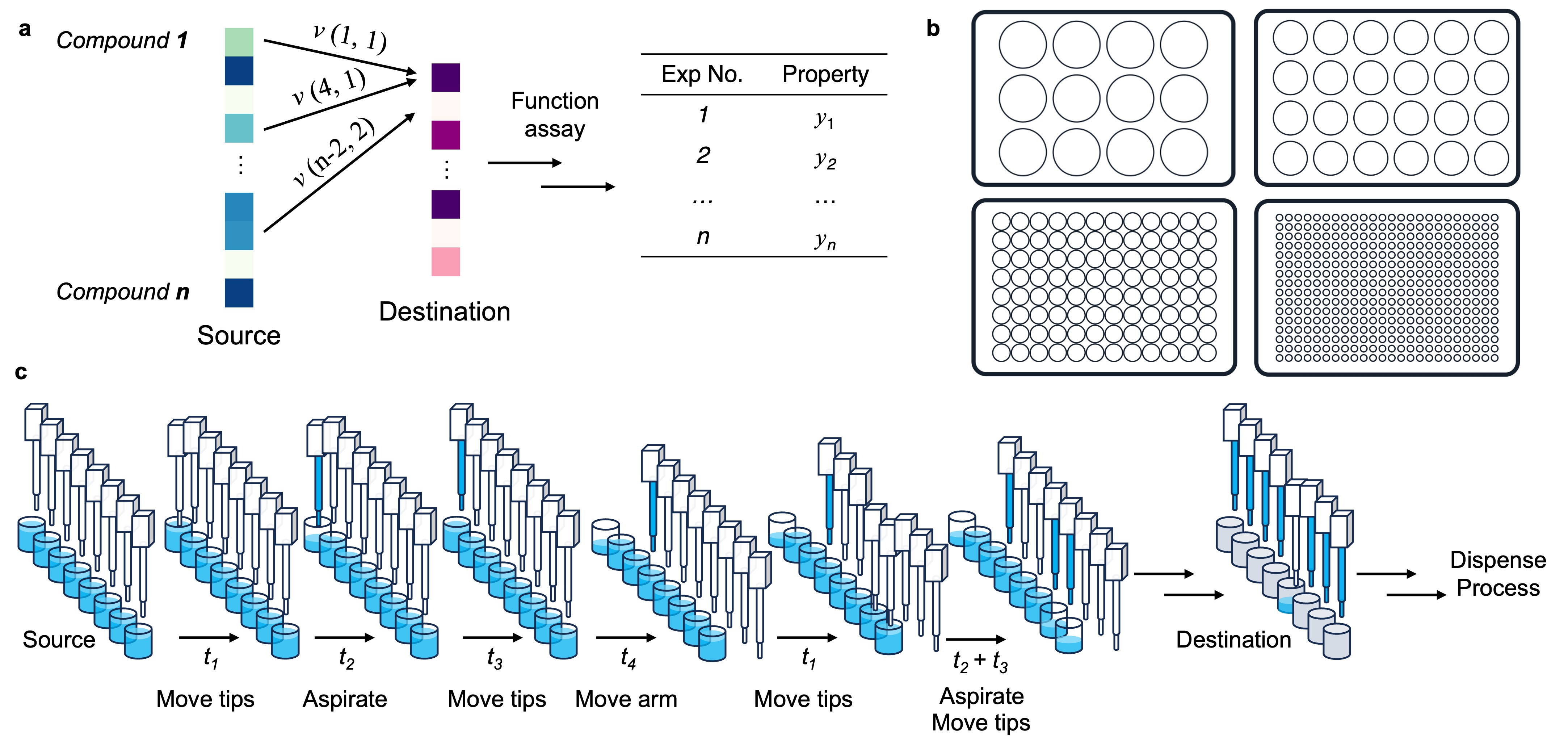}
 \caption{Problem description. (a) Schematic representation of a typical multi-channel pipetting task, showing how different compounds and volumes are allocated from source wells to destination wells for downstream functional assays.
(b) Commonly used labware formats, including 12-, 24-, 96-, 384-well plates.
(c) Step-by-step illustration of the liquid handling process. The time required can vary significantly based on the liquid characteristics and volumes.}
 \label{problem_des}
\end{figure*}

Optimizing (reducing) execution time would lead to substantial improvements in throughput and efficiency.
Among the various liquid handling platform, the 8-channel pipette stands out as one of the most widely used configurations. Of the available 8-channel pipette configurations, individually addressable pipettes offer superior flexibility, making them well-suited for combinatorial formulation screening. Such pipette can be found in Tecan, Hamilton, Beckman, Revvity, and other liquid handling platforms.
Experimental protocols defining the precise sequence and order of liquid transfer operations are typically defined by a user without explicit optimization of execution time. 
Despite the ubiquity of liquid handling operations, to the best of our knowledge, no existing method in the literature offers an approach to systematically optimizing execution time of the liquid handling tasks. And despite its seeming simplicity, this combinatorial pipette scheduling problem is non-trivial and offers substantial room for efficiency gains.

\begin{figure*}
 \centering
 \includegraphics[height=15cm]{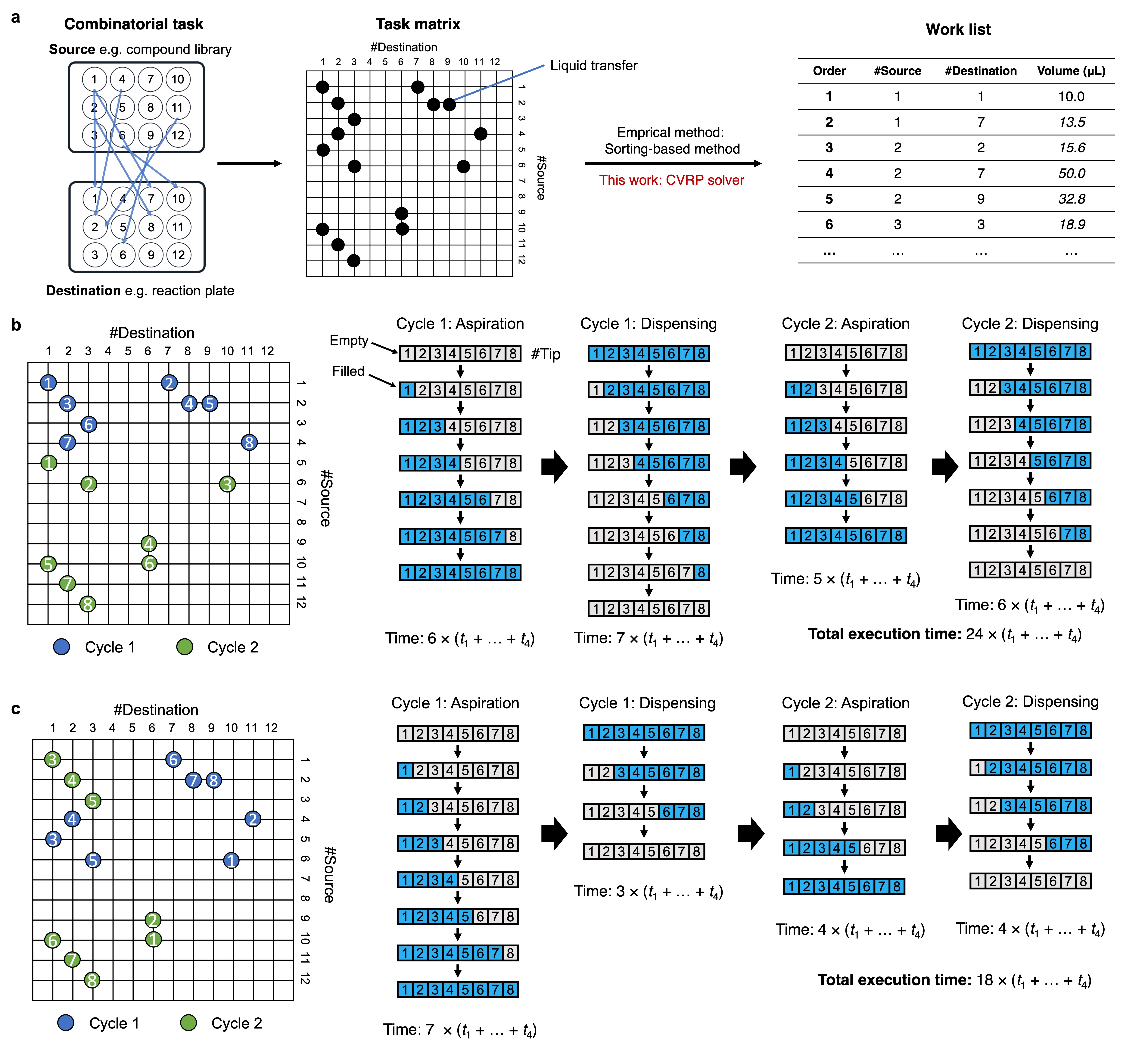}
 \caption{
Problem formulation. (a) The workflow from the defined combinatorial task between two 12-well plates to executable work list. In this work, the position is numbered follow column-major indexing. For example, A1 is well 1, B1 is well 2 \textit{etc.} After specifying the labwares and the arrangement of the reagents and experiments, the task is represented as a task matrix, non-zero entries denote individual liquid transfer—each defined by (position on source plate, position on destination plate, volume). Different scheduling methods can be used to generate the work list in different orders. We apply a CVRP-based solver to derive a work list that minimizes the number of tip movements, thereby reducing overall execution time. (b-c) Demonstration of the influence of the work list order on the total execution time of the same task involving 16 liquid transfers with (b) a less efficient and (c) a more efficient execution sequence. The task requires 2 cycles of aspiration-dispensing with 8-channel pipette. The numbers in cycle are the orders of the liquid transfer in the indicated cycle. Each cube represents a pipette tip, labeled 1–8. Blue indicates a tip filled with liquid, while gray indicates an empty tip. $t_1$ to $t_4$ represents the time required for tip lowering, aspirating/dispensing, withdrawing, and moving to the next location. The time for aspiration and dispensing are taken to be equal in this illustration.}
 \label{vrp_demo}
\end{figure*}

Herein, we propose an optimization strategy to systematically reduce the execution time of liquid handling tasks on 8-channel systems with individually controllable tips. Our key contributions include: (1) defining a function that serves as a robust proxy for the execution time, and (2) formulating the scheduling challenge as a Capacitated Vehicle Routing Problem (CVRP), which enables the use of heuristic solvers traditionally applied in logistics and transportation planning.
This approach significantly improved the efficiency of task planning and execution, resulting in a up to 37\% performance improvement compared to the baseline sorting method. The results underscore the substantial potential for optimizing the operation of existing liquid handling platforms without changing the hardware configuration, paving the way for more efficient high-throughput experimentation and better utilization of the growing repertoire of autonomous laboratories.

\section{Methods}
\subsection*{Problem description.~~}
The combinatorial liquid handling task involves transferring varying volumes of multiple compounds from a set of sources to designated destinations, following a predefined experimental design (\Cref{problem_des}a). We focus on executing these tasks using liquid handling system with 8 individually controllable channels, a widely adopted configuration in laboratory automation. Standard Society for Biomolecular Screening (SBS)-format well plates, including 12-well, 24-well, 96-well, and 384-well plates, are commonly used for storing source compounds and receiving reagents (\Cref{problem_des}b). These formats differ in layout and spacing, requiring different pipetting strategy. For instance, a 96-well plate allows all eight pipette tips to simultaneously access a single column, whereas a 12-well plate accommodates only three tips per column due to its 3 (rows) $\times$4 (columns) layout. While a 384-well plate supports eight tips per column, the narrower spacing restricts tip placement to every other well to avoid physical collisions. These geometric constraints must be considered during pipette scheduling to ensure accuracy, efficiency, and compatibility with the selected labware format.

The liquid handling operation consists of a sequence of discrete actions (\Cref{problem_des}c), including lowering the tip into the liquid ($t_1$), aspirating or dispensing ($t_2$), withdrawing the tip ($t_3$), and moving the arm to the next location ($t_4$). 
Given $n$ ($n\leq 8$) available tips, the mainstream liquid handling platforms typically perform the liquid transfer based on a work list consisting of (source, destination, volume) entries, executing them in $n$-by-$n$ batches (\Cref{vrp_demo}a), followed by washing (for fixed tips) or tip replacement (for disposable tips) after completion of the dispensing operation. Each step incurs a time cost, which can vary significantly based on the layout of the bench, volume, liquid viscosity and required accuracy. As mentioned, the arm movement time ($t_4$) is typically relatively smaller compared to others; tip lowering ($t_1$) and withdrawal ($t_3$) times are usually similar to each other in duration; and aspiration or dispensing times ($t_2$) depend on the transfer volume in addition to material properties. For instance, viscous liquids typically demand slow aspiration, dispensing, and withdraw time to maintain volume  precision\cite{thurow2022devices}. 
Depending on the order of the work list, the same liquid handling task could have significantly different execution time (\Cref{vrp_demo}b,c). Maximizing tip movement parallelization reduces the number of lower–aspirate/dispense–withdraw–move cycles required, thereby enabling a more efficient liquid transfer process (\Cref{vrp_demo}c). One effective strategy is to maximize the number of next-tip tasks which use adjacent tips for aspiration or dispensing, through optimizing the order of the work list. In this way, we can minimize arm movements while filling or emptying all available channels. %

\begin{figure}[t]
\centering
  \includegraphics[height=6cm]{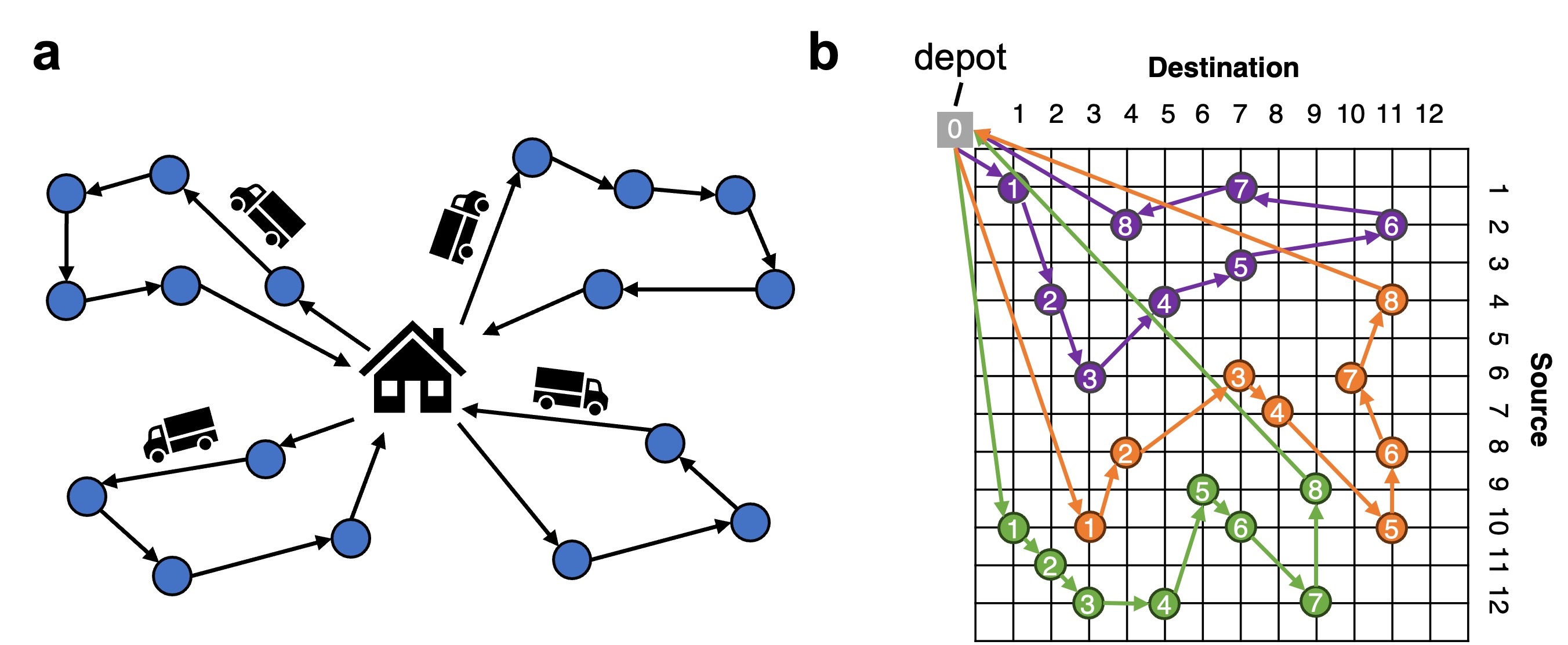}
  \caption{The analogy between capacitated vehicle routing problem (CVRP) and the scheduling challenge of liquid handling. (a) CVRP is a variant of the vehicle routing problem in which each vehicle has a limited carrying capacity. The objective is to determine the most cost-efficient set of routes that service all locations without exceeding the capacity constraints of any vehicle. (b) In liquid handling, an 8-channel pipette (analogous to a vehicle with 8-unit capacity) must perform multiple source-to-destination liquid transfers (analogous to locations). Each aspiration-dispensing cycle corresponds to a delivery route, and the goal is to minimize the total execution time. This analogy enables the use of CVRP solvers to globally optimize pipette scheduling and reduce execution time.}
  \label{cvrp_analogy}
\end{figure}

This scheduling challenge bears a strong resemblance to the Capacitated Vehicle Routing Problem (CVRP), a classical combinatorial optimization problem in operations research (\Cref{cvrp_analogy}a). In CVRP, a fleet of vehicles need to determine the most efficient routes to deliver goods to a set of locations, starting and ending at a central depot, while minimizing total travel cost and satisfying constraints such as vehicle capacity. Drawing an analogy to pipette scheduling, each (source, destination) pair can be viewed as a location to be visited (\Cref{cvrp_analogy}b), and the 8-channel pipette functions as a vehicle with a capacity of 8 deliveries per cycle. The “distance” between locations is defined by their relative positions on both the source and destination plates, as determined by the physical geometry of the well plate. Wells aligned in the same column and adjacent rows are considered closer and more efficient to access within a single operation. Importantly, this spatial relationship is directional. For example, within a column of a 96-well plate, row 3 is close to row 4 but not to row 2. This avoids misaligned assignments, such as tip 3 aspirating from row 2 and tip 2 from row 3, which would otherwise lead to unnecessary arm movements. By framing the problem in this way, we can apply CVRP solvers to minimize the total computed (estimated) execution time. 

\subsection*{Mathematical formulation of the scheduling challenge.~~}
We denote the action of aspirating from well $a$ and dispensing to well $b$ as a single job of the scheduling task.
Jobs are encoded as non-zero entries in the task matrix (\Cref{vrp_demo}a) $\mathbf{T} \in \mathbb{R}^{n^{src}\times n^{dst}}$, where $t_{a,b}$ denotes the volume, $n^{src}$ is the number of wells in the source plate, and $n^{dst}$ is the number of wells in the destination plate. 
Solving the scheduling task is equivalent to finding the optimal sequence of executing all jobs that minimizes the total time cost required to finish a pipetting task.

We first define the pairwise distance of aspirating or dispensing two wells consecutively. We define a unit action as moving tip, aspirating/dispensing, moving tip again, and moving arm ($t_1+t_2+t_3+t_4$, \Cref{problem_des}c). In this work, we ignore the impact of different distances when moving arms, as we noticed that they are not dominating factors in the time cost for commercially available liquid handlers we have used. 
For a plate with $n$ wells (e.g., $n=96$), we define the following pairwise distance matrix $\mathbf{D} \in \{0,1\}^{n\times n}$:
\begin{equation}
    d_{a,b} =\left\{\begin{array}{cl}
       0  & \text{if } a<b \text{ and } a,b \text{ are adjacent wells in the same row;}\\
       1  & \text{otherwise.}
    \end{array} \right.
\end{equation}
Wells next to each other can be aspirated or dispensed at the same time, meaning that when these two jobs are ranked consecutively in the work list, there is no extra cost for the liquid handler as they will in practice be executed simultaneously. Due to differences in well spacing, adjacency is defined differently for higher-density plates: for a 384-well plate, adjacent wells correspond to every other well in the row; for a 1536-well plate, adjacency occurs every four wells. If not, another unit operation is needed to finish these two jobs. We compute $\mathbf{D}^{src}$ and $\mathbf{D}^{dst}$ for the source plate and the destination plate, respectively. 

Recall that dispensing well $a$ in the source plate to well $b$ in the destination plate is defined as a job. Our next step is to construct a job-level distance matrix. Assuming we have $m$ jobs, we define $\mathbf{S}\in \{0,1\}^{m\times n^{src}}$ and $\mathbf{E}\in \{0,1\}^{m\times n^{dst}}$ as the incidence matrices of $\mathbf{T}$, where  
$s_{i,a} = 1, e_{i,b} = 1$
if task $i$ is to aspirate from well $a$ in the source and dispense to well $b$ in the destination. With $\mathbf{S}$ and $\mathbf{E}$, we are able to transform the pairwise distances for each well to the following job-level distance matrices,
\begin{equation}
    \bar{\mathbf{D}}^{{src}^{\prime}} = \mathbf{S}\mathbf{D}^{src}\mathbf{S}^\top, \quad \bar{\mathbf{D}}^{{dst}^{\prime}} = \mathbf{E}\mathbf{D}^{dst}\mathbf{E}^\top.
    \label{eq:dprime}
\end{equation}

We define a new matrix ${\mathbf{D}}^{\prime}\in\mathbb{R}^{(m+1)\times(m+1)}$, where index 0 corresponds to a dummy job, as 
${\mathbf{D}}^{\prime}={\mathbf{D}}^{src^\prime} + {\mathbf{D}}^{dst^\prime}$, with 
\begin{align}
    {d}^{src^\prime}_{i,j} &= \left\{\begin{array}{ll}
        t_{1,3,4}^{src}+ \frac{v_{j}}{q^{src}} & \text{if } \bar{d}^{src^\prime}_{i,j}=1 \text{ and } i\geq0,j\geq 1,\\
        \max\left(0, \frac{v_{j} - v_{i}}{q^{src}}\right) & \text{otherwise,}
    \end{array}\right. \\
    {d}^{dst^\prime}_{i,j} &= \left\{\begin{array}{ll}
        t_{1,3,4}^{dst}+ \frac{v_{j}}{q^{dst}} & \text{if } \bar{d}^{dst^\prime}_{i,j}=1 \text{ and } i\geq0,j\geq 1,\\
        \max\left(0, \frac{v_{j} - v_{i}}{q^{dst}}\right) & \text{otherwise,}
    \end{array}\right.
\end{align}
where $t^{src}_{1,3,4}$ (s) is the sum of $t_1,t_3,t_4$ for aspirating at the source plate, which is (approximately) viewed as a constant. $v_j$ ($\mu$L) is the volume for job $j$, and $q^{src}$ ($\mu$L/s) is the speed of aspiration, therefore $\frac{v_j}{q^{src}}$ is the aspiration time ($t_2$) for job $j$. The same definitions are applied to dispensing operations. The dummy job has $v_0 = 0$.
The pipette scheduling problem can then be formulated as,
\begin{subequations}
\label{eq:cvrp_form}
\begin{align}
    \min_{\mathbf{X}} \ &\sum_{i = 0}^m\sum_{j=0}^m \sum_{k=1}^{K} d^\prime_{i,j} \cdot x_{i,j,k} \label{eq:cvrp_obj}\\
    s.t.\quad &  \sum_{i=0}^m x_{i,j,k} = \sum_{i=0}^m x_{j,i,k}, \label{eq:leaving_entering}\\ 
    & \sum_{k=1}^K \sum_{i=0}^m x_{i,j,k} = 1 \quad \forall j \{1...m\}, \label{eq:node_once}\\ 
    & \sum_{j=1}^m x_{0,j,k} = 1\quad \forall k\in\{1...K\}, \label{eq:leaves_depot}\\ 
    & \sum_{i=0}^m \sum_{j=1}^m x_{i,j,k} \leq 8 \quad \forall k\in\{1...K\}, \label{eq:capacity}\\ 
    & \mathbf{X}\in\{0,1\}^{(m+1)\times(m+1)\times K}, \\
    & x_{i,i,k} = 0 \quad \forall i \in \{0 ...m\}, k\in\{1...K\}.
\end{align}
\end{subequations}
$\mathbf{X}$ is the aspirating and dispensing plan, where $x_{i,j,k}=1$ means job $i$ is followed by job $j$ at cycle $k$. $K=\lceil \frac{m}{8}\rceil$ denotes the number of cycles needed to dispense all jobs, because the liquid handler can aspirate or dispense at most 8 wells at the same time; one could easily generalize to non-conventional liquid handlers by changing this number. 
\Cref{eq:cvrp_obj} is the \textit{computed execution time} of the pipette task.  
Constraint (\ref{eq:leaving_entering}) means the number of times leaving a job should be the same as the number of times entering a job, and constraint (\ref{eq:node_once}) ensures each job is completed exactly once. Constraint (\ref{eq:leaves_depot}) means that each cycle should leave the dummy job, which helps enforce constraint (\ref{eq:capacity}) that the capacity of each cycle is 8 jobs. 

\subsection*{Solver implementation.~~}

The formulation in \Cref{eq:cvrp_form} is exactly the same as CVRP, where the dummy job ($i=0,j=0$) is treated as the shared vehicle depot, $K$ cycles are equivalent to $K$ vehicles, each vehicle has a capacity of 8, and the distance matrix $\mathbf{D}^\prime$ is interpreted as the pairwise routing distance. To this end, we can tackle pipette scheduling with off-the-shelf CVRP solvers.

The pipette scheduler is developed with the CVRP solver implemented in Google OR-Tools~\cite{ortools_routing}. All the computation in this work was performed on a laptop (MacBook Pro with M3 Pro, 18GB RAM). We compute $\mathbf{D}^\prime$ from the plate layout and liquid handling parameters ($t^{src}_{1,3,4}$, $q^{src}$, $t^{dst}_{1,3,4}$, and $q^{dst}$) and job based on \Cref{eq:dprime}, and pass $\mathbf{D}^\prime$ to the CVRP solver as the distance matrix. Each job is viewed as a location to visit in CVRP, and the dummy node with all-zero distances to all other locations is defined as the depot (i.e., starting location) in CVRP. After getting the routing result from the solver, we translate it into the corresponding pipette work list under a format that the liquid handler control software can parse (\Cref{vrp_demo}a). Unless otherwise specified, the $t^{src}_{1,3,4}$ and $t^{dst}_{1,3,4}$ were set to 1, and $q^{src}$ and $q^{dst}$ were set to 100 in the subsequent results.

\subsection*{Baseline methods.~~}
We evaluate the performance of our CVRP-based scheduling approach by comparing it with heuristic baseline methods. The first method, named long-axis prioritized (LAP) method, is a parallelization-driven strategy that attempts to maximize the number of simultaneous transfers on the plate (source or destination plate) with more number of wells by iteratively sampling the jobs on the axis that belongs to the larger plate until all of the jobs are sampled. This method could guarantee partial parallelization on at least one of the labware components. The second, named greedy, is a greedy heuristic that randomly selects the closest (source, destination) jobs on distance matrix (${\mathbf{D}}^\prime$) to iteratively pick the nearest unassigned pair. Additionally, we included a control method, named row-major sorting, where jobs are executed in the order returned by \texttt{np.argwhere(T)}, which corresponds to row-major order—i.e., traversing the task matrix from top to bottom and left to right. 

\subsection*{Random task generation.~~}

To generate synthetic pipetting tasks for benchmarking, we implemented a custom random sampling procedure. Given specified source and destination plates, we initialized a two-dimensional matrix of zeros with shape ($n^{src}$, $n^{dst}$), where $n$ is the number of wells of the labware. We randomly selected a defined number of unique positions in the matrix to assign non-zero values, corresponding to the liquid transfers. The number of non-zero elements reflects the total number of transfers in the task. Each selected position was assigned to a random volume sampled uniformly between 1 and 100 as a representation of the volume to be transfered. 

\subsection*{Simulation of the tasks.~~}

Liquid handling task execution was simulated using EvoSim software (version 2.8.0.0, Tecan) and the \textit{simulated execution time} was calculated by subtracting the start time from the end time of each simulated run. The software accepts a .csv input file containing a list of pipetting instructions, each defined by a (source position, destination position, volume) triplet. These instructions are executed on a virtual worktable in simulation mode with 3D rendering (Figure S1).

Simulations were performed in normal speed mode, which could reflect the real-world execution time. After each aspiration–dispensing cycle, an additional washing step was included. The detailed configuration of the worktable layout is shown in Figure S1, and all operational parameters used in the simulations are provided in Table S1.

\section{Results and discussion}

\subsection*{The computed execution time is a robust proxy of the execution time.~~}

\begin{figure}[h]
\centering
  \includegraphics[height=10cm]{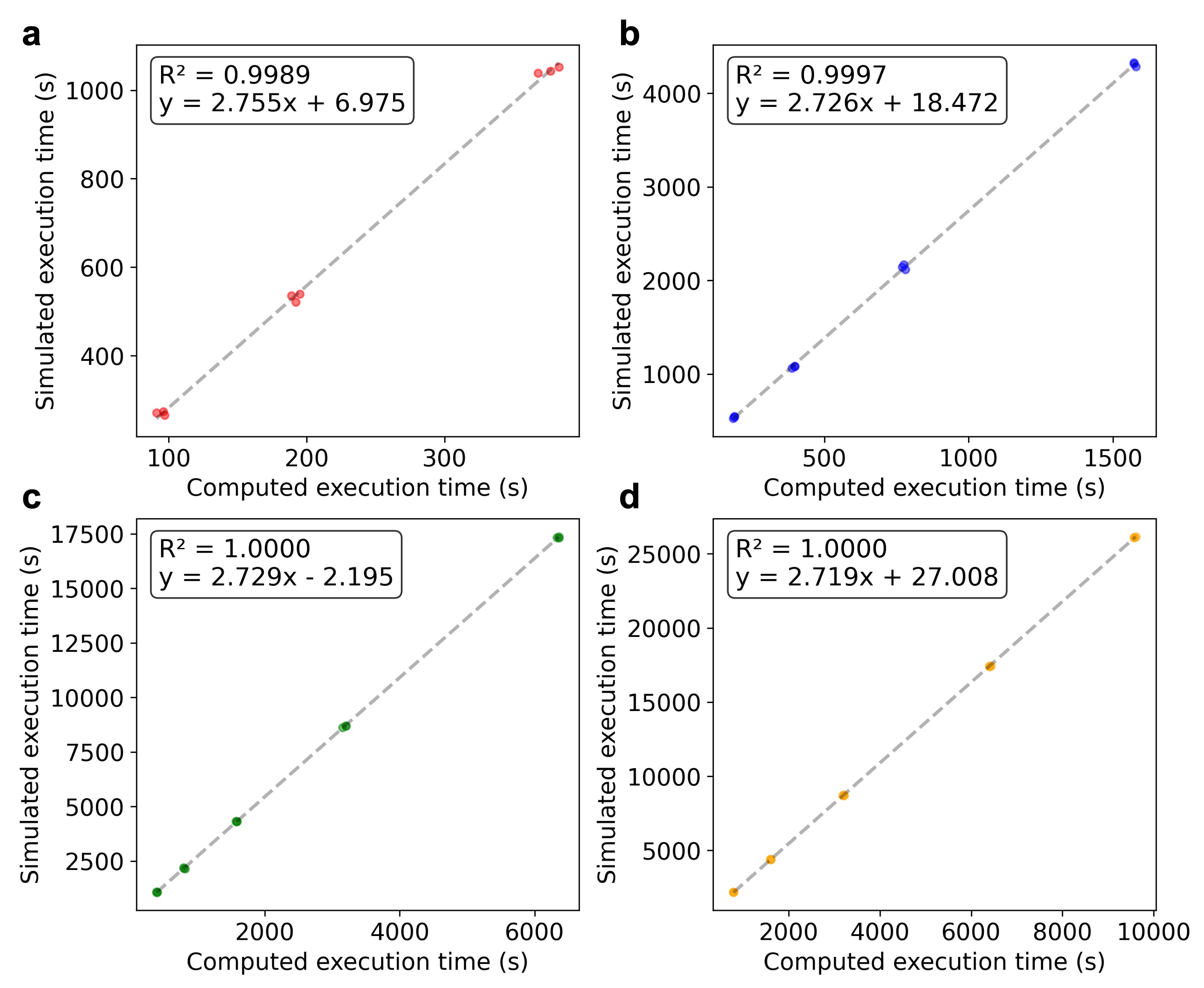}
  \caption{The correlation between the computed execution time and the simulated execution time of liquid handling tasks between (a) a 12-well plate and a 12-well plate with 25, 50 and 100 transfers, (b) a 24-well plate and a 24-well plate with 50, 100, 200 and 400 transfers, (c) a 96-well plate and a 96-well plate with 100, 200, 400, 800 and 1600 transfers, and (d) a 384-well plate and a 384-well plate with 200, 400, 800, 1600 and 2400 transfers. For each labware and number of transfers, 3 random task matrices were generated and evaluated. Correlation coefficients and $R^2$ are shown for each. }
  \label{proxy}
\end{figure}

To assess whether the computed execution time could serve as a reliable proxy for actual execution time, we performed a series of simulations of randomly generated tasks (see Methods) with a specified number of liquid transfers between same type of labware in \Cref{problem_des}b. For each task, the work list was constructed by randomly ordering the transfer operations. This setup provides an unbiased framework for systematic evaluation across a broad range of task arrangement, ordering, and labware formats. Importantly, this approach ensures that performance comparisons are not influenced by the structure or assumptions of any specific experimental protocol, thereby enabling generalizable insights into algorithmic effectiveness. As shown in \Cref{proxy}, a strong correlation was observed between the computed execution time and the simulated execution time. 
The results confirmed that our definition of the computed execution time can be used as a proxy for estimating execution time during pipette scheduling; minimization of the former should lead to minimization of the latter.

\subsection*{The CVRP-based method consistently outperforms other methods in terms of execution time of the proposed pipetting strategy.~~}

\begin{figure*}[h]
 \centering
 \includegraphics[height=9cm]{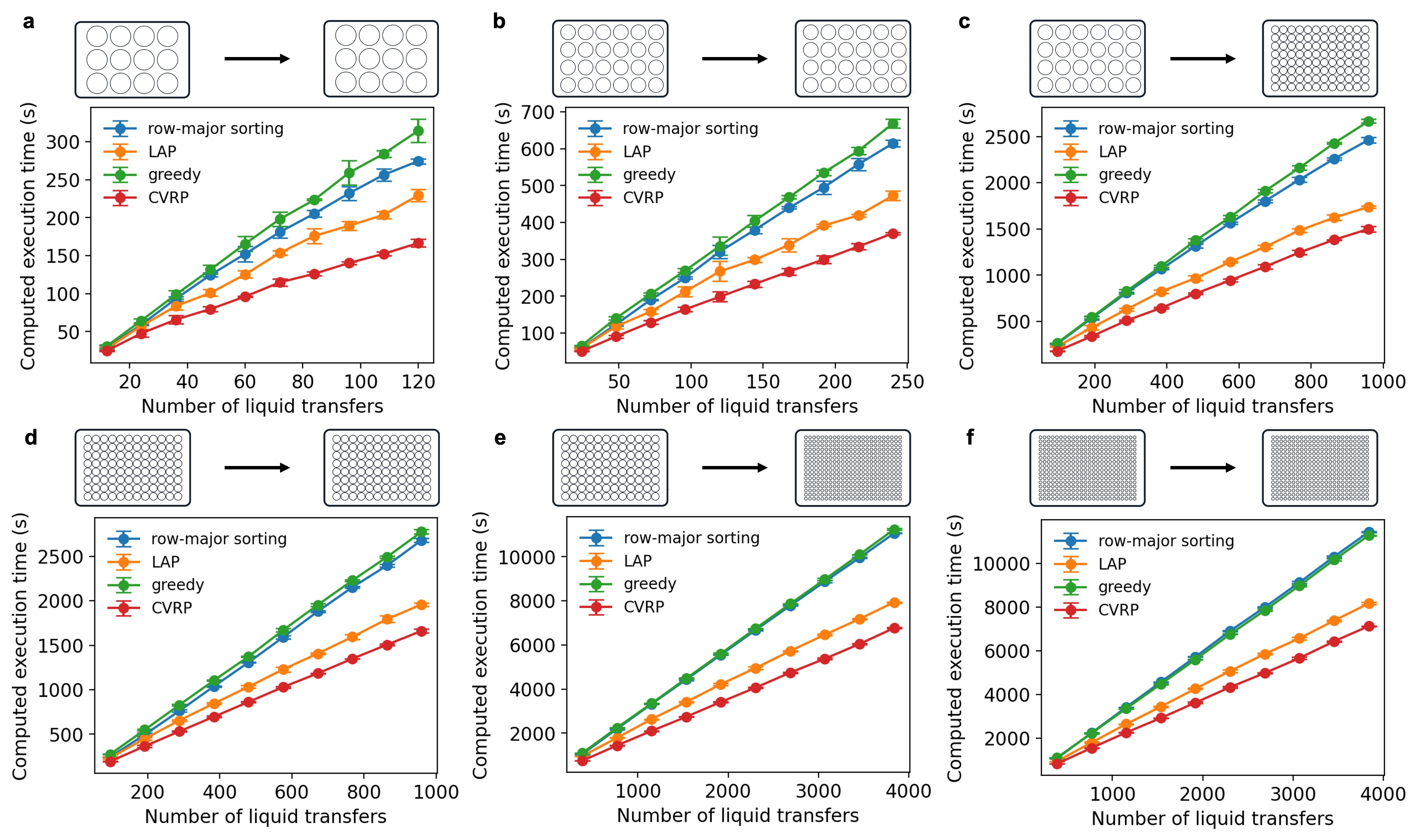}
 \caption{Comparison of the performance of different scheduling methods on randomly generated tasks (a) from 12-well plate to 12-well plate, (b) from 24-well plate to 24-well plate, (c) from 24-well plate to 96-well plate, (d) from 96-well plate to 96-well plate, (e) from 96-well plate to 384-well plate and (f) from 96-well plate to 96-well plate. For a given number of liquid transfers and labware types, 3 random tasks were generated and scheduled using different methods. The solving time for CVRP is 20 seconds. For each labware combination, we evaluated 10 different task sizes, with the number of liquid transfers set to multiples (from 1× to 10×) of the number of wells in the destination plate. Results are shown as mean ± standard deviation. The results of the remaining combinations are provided in Figure S2.}
 \label{optimization_results}
\end{figure*}
We observed that the CVRP-based scheduling method consistently outperforms baseline methods in minimizing execution time for randomly generated tasks (\Cref{optimization_results} and Figure S2. 
The performance is robust across various labware formats and remains effective for tasks involving up to approximately 4,000 liquid transfers. On average, the CVRP-based method achieved a 37\% reduction in execution time compared to the row-major sorting method. While the LAP method exhibited near-optimal performance in certain labware combinations (\Cref{optimization_results}c–e), its overall performance was inconsistent. For instance, in lower-density formats such as 12- and 24-well plates, the improvement over row-major sorting methods was less pronounced compared to higher-density formats like 96- and 384-well plates. In contrast, the CVRP-based approach consistently delivered performance gains across all formats, reducing execution time by an average of 15\% relative to the LAP method. These findings underscore the robustness and generalizability of the CVRP formulation, particularly in scenarios where baseline heuristics may fail to provide consistent improvements.

\begin{figure}[h]
\centering
  \includegraphics[height=6cm]{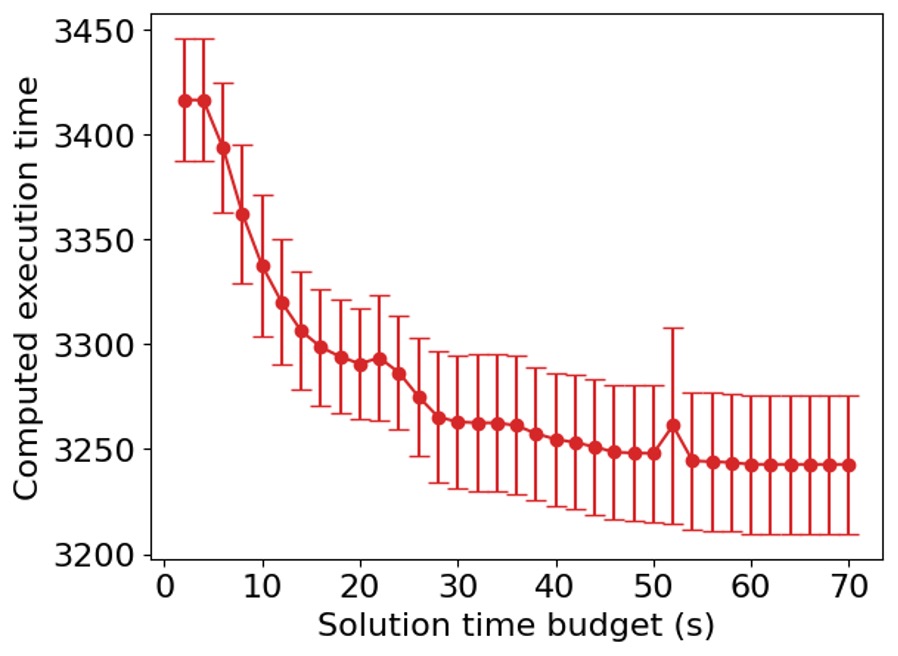}
  \caption{Relationship between solver performance and solution time budget. 6 random tasks, each with 2000 liquid transfers, were generated and solved with different solution time budget from 2 to 70 seconds, in 2-second increments. Results are presented as mean ± standard deviation.}
  \label{time_performance}
\end{figure}

We next investigated how the solution time allocated to the solver affects optimization performance. CVRP is an NP-hard problem for which it is impractical to find the optimal solution in polynomial time; hence, we resort to the approximate solver in OR-Tools. The solver requires a minimum amount of time to produce a feasible solution; if insufficient time is allocated, the program may fail to return a result. To further explore the relationship between the solution time and optimization effect, we evaluated solver performance on tasks involving 2,000 liquid transfers from a 96-well plate to another 96-well plate. As shown in \Cref{time_performance}, increasing the allotted solution time consistently improved performance until a plateau was reached starting at around 40 CPU seconds. 

The proposed method can be readily generalized to high-density labware such as 1536-well plates. With a solution time of 120 CPU seconds, e successfully optimized pipetting tasks involving up to approximately 14,000 liquid transfers (Figure S3). For the most complex task evaluated, the method reduced the computed execution time to 25565 compared to the 29042 of the LAP method, corresponding to an execution time reduction of 158 minutes relative to the LAP scheduling strategy. These results demonstrate the scalability of the approach and its potential to deliver substantial time savings in large-scale, high-throughput liquid handling operations.

\subsection*{Optimizing the schedule of real-world tasks leads to tangible improvements in efficiency and throughput.}

\begin{figure}[h]
\centering
  \includegraphics[height=11cm]{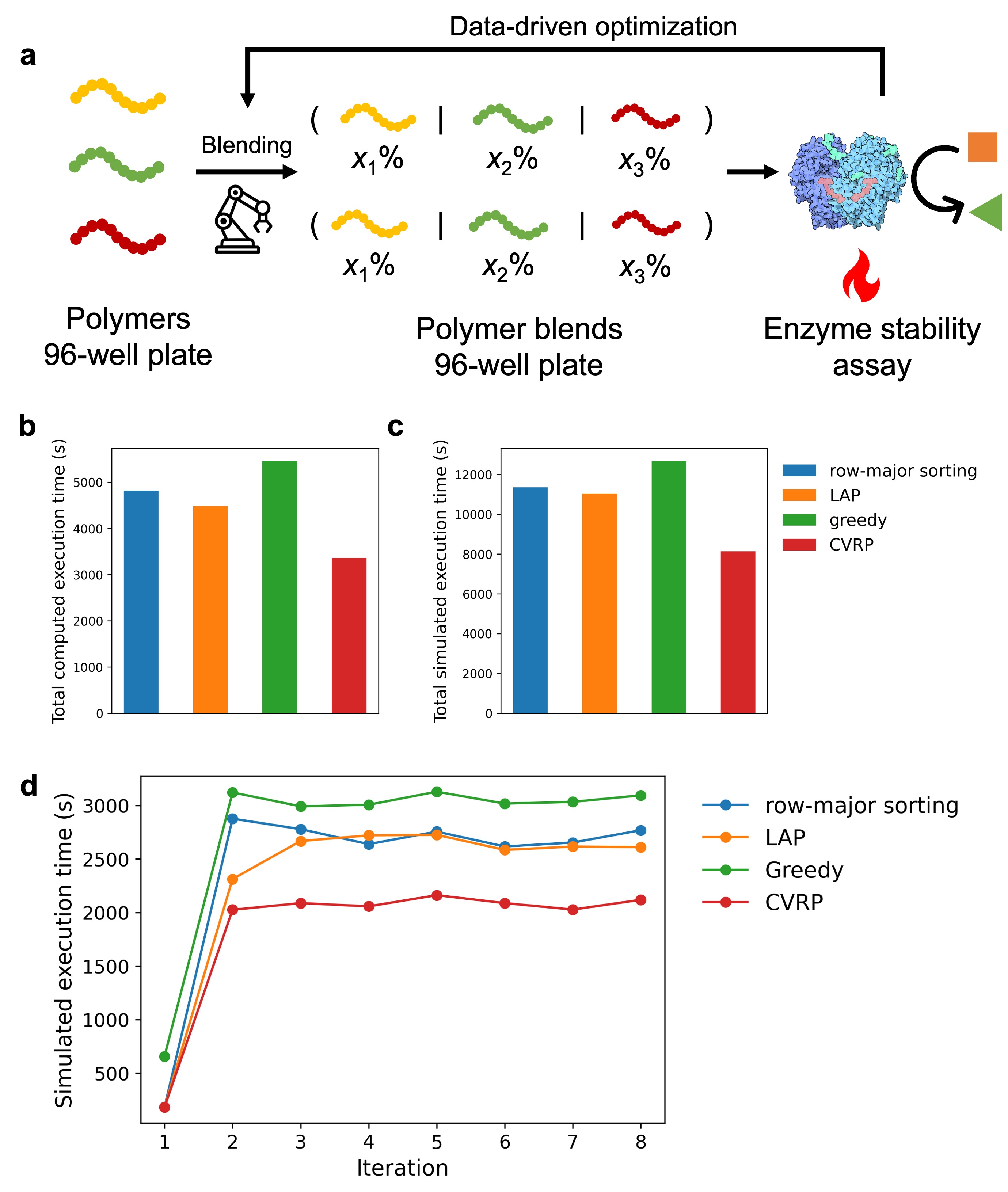}
  \caption{Optimization of the execution time of a real-world task.
  (a) Schematic demonstration of the an autonomous polymer blend optimization campaign for enzyme stabilization. The polymer stock solutions were stored in the source 96-well plate and were blended at the destination 96-well plate. The enzyme stabilization assay is then performed in the destination well plate. The new experiment for the next round is proposed by a genetic algorithm. (b) Total computed execution time and (c) total simulated execution time of different methods.
(d) Simulated execution time of each iteration across campaign using different methods. Solving time for CVRP is 20 seconds. }
  \label{real_world}
\end{figure}

We then demonstrated the performance of our method on a real-world task derived from a previously developed automated experimental platform for the discovery of random heteropolymer blends for enzyme stabilization\cite{wu2024autonomous}. This workflow (\Cref{real_world}a) involves high-dimensional combinatorial liquid transfer operations to blend polymer stock solutions from one 96-well plate to another. Through this self-driving platform, we successfully identified polymer blends that outperformed their individual constituents and can stabilize the glucose oxidase under 70°C for 30 minutes. However, the exploration capacity of the autonomous platform is constrained by the time required to execute the blending process. As the number of components in each blend increases, the associated liquid handling time grows significantly—often exceeding the practical limits imposed by the shelf life of sensitive reagents such as enzymes. This limitation was a key factor in our decision to restrict the number of blend components to 4. Improving scheduling efficiency could enable more experiments within the same time frame or allow exploration of a larger design space without compromising reagent stability.

This real-world task is different from the randomly generated tasks. The blending composition of the polymer stock solutions to be added to the destination wells are proposed by an optimization algorithm based on the outcomes of previous iterations. As the experiment progresses, certain source wells become increasingly favored or disfavored, resulting in a non-uniform distribution of liquid transfers on the task matrices. Additionally, the layout of the destination well plate must conform to specific rules to accommodate control experiments, further complicating scheduling (Figure S4).

We optimized the polymer blending process of one real experimental campaign from this work. In the original workflow, work lists were generated using the row-major sorting method and executed on a Tecan Evo 200 liquid handling platform. The campaign started from a control experiment whose task matrix is a diagonal matrix with first 8 rows empty for the control experiments (Figure S4), followed by a round of pure random exploration. Subsequent experiments were generated adaptively by a genetic algorithm based on prior results. For the CVRP-based method, we allocated 20 seconds of solving time per iteration, totaling approximately 3 minutes for the entire campaign. As shown in \Cref{real_world}b, the CVRP-based approach significantly outperformed all other methods in reducing total computed execution time. Notably, it achieved a 25\% reduction compared to the LAP method—substantially greater than the 15\% reduction observed in purely random tasks. Execution time simulations further supported this finding (\Cref{real_world}c), the CVRP-based method achieved a total simulated execution time of 246 minutes, compared to 307 minutes with the LAP method and 321 minutes with the row-major sorting method, representing time savings of 61 minutes and 75 minutes, respectively (throughput improvements of 25\% and 30\%). This improvement is attributed to the reduced effectiveness of the LAP method in handling non-random task starting from the third iteration (\Cref{real_world}d). The results underscore the robustness and effectiveness of the CVRP-based optimization, particularly in dynamic, data-driven workflows where traditional heuristics fail to perform consistently.

To evaluate the generalizability of our strategy across different liquid handling platforms, we tested the task of iteration 3 in \Cref{real_world}d---the iteration where we observe the proposals from different scheduling methods to diverge greatly in simulated execution time---on a JANUS G3 automated liquid handling workstation (Revvit) with different aspiration and dispensing speeds. The CVRP-based method outperformed all other methods in all the speed combinations. At an aspiration speed of 100 $\mu$L/s and dispensing speed of 25 $\mu$L/s, it achieved an execution time of 36 minutes compared to the LAP method's 45 minutes (Figure S5). This result further demonstrates the versatility and platform-independence of our approach, underscoring its potential to improve efficiency across a wide range of automated systems.

\section{Conclusions}
We have demonstrated how the execution time of 8-channel liquid handling tasks can be effectively optimized as a Capacitated Vehicle Routing Problem, or CVRP. We achieved substantial reductions in execution time in both simulated and experimental settings through a lightweight optimization step, which enables gains in throughput. 
Further improvements might be realized by incorporating layout-aware destination assignment strategies during experimental design to further reduce execution overhead. As laboratory automation continues to play a pivotal role in accelerating scientific discovery, our method provides a practical, scalable, and generalizable solution for improving throughput and efficiency without hardware modification. It can be readily integrated into formulation optimization platforms and other high-throughput experimental workflows involving combinatorial screening of chemical or biological systems.
\section*{Conflicts of interest}

There are no conflicts to declare.
\section*{Author contributions}

G.W. and R.W. conceived the project. G.W. identified the problem, developed the code, and performed the simulations. R.W. contributed to the codebase and formulated the problem mathematically. C.W.C. supervised the project. All authors contributed to writing the manuscript.

\section*{Data availability}

The code is available at: \url{https://github.com/wuRoy/CVRP_pipette_scheduling.git}

\begin{acknowledgement}

G.W. acknowledges the discussion with Dirk Lauinger and David Willshaw for the problem formulation. We thank Haisen Zhou (Peking University) for testing the script on JANUS G3 automated liquid handling workstations. Research reported in this publication was supported by NIGMS of the National
Institutes of Health under award number R21GM141616. The content is solely the responsibility of the
authors and does not necessarily represent the official views of the National Institutes of Health.

\end{acknowledgement}

\begin{suppinfo}

\end{suppinfo}

\bibliography{achemso-demo}

\end{document}


\maketitle

\tableofcontents
\newpage

\section{Simulation and Operation Parameters}
All simulations were conducted using EvoSim. Key operational parameters are summarized in Table~\ref{sim_params}.

\begin{table}[h]
\centering
\caption{Key liquid handling simulation parameters in EvoSim.}
\label{sim_params}
\begin{tabular}{ll}
\hline
\textbf{Parameter}           & \textbf{Value}     \\
\hline
Aspiration speed             & 100~$\mu$L/s       \\
Aspiration delay             & 500~ms             \\
Dispensing speed             & 100~$\mu$L/s       \\
Dispensing delay             & 500~ms             \\
Tip retraction speed             & 20 mm/s             \\
\hline
\end{tabular}
\end{table}

\begin{table}[h]
\centering
\caption{Key liquid handling parameters of JANUS liquid handling workstation.}
\label{janus_params}
\begin{tabular}{ll}
\hline
\textbf{Parameter}           & \textbf{Value}     \\
\hline
Aspiration speed             & As specified      \\
Aspiration delay             & 500~ms             \\
Dispensing speed             & As specified     \\
Dispensing delay             & 500~ms             \\
Scan in speed             & 150 mm/s             \\
Scan out speed             & 150 mm/s             \\
Retract from liquid speed             & 100 mm/s             \\
Retract from liquid hight             & 10 mm             \\
\hline
\end{tabular}
\end{table}

\section{Additional Figures}

\begin{figure}[H]
\centering
\includegraphics[width=0.85\textwidth]{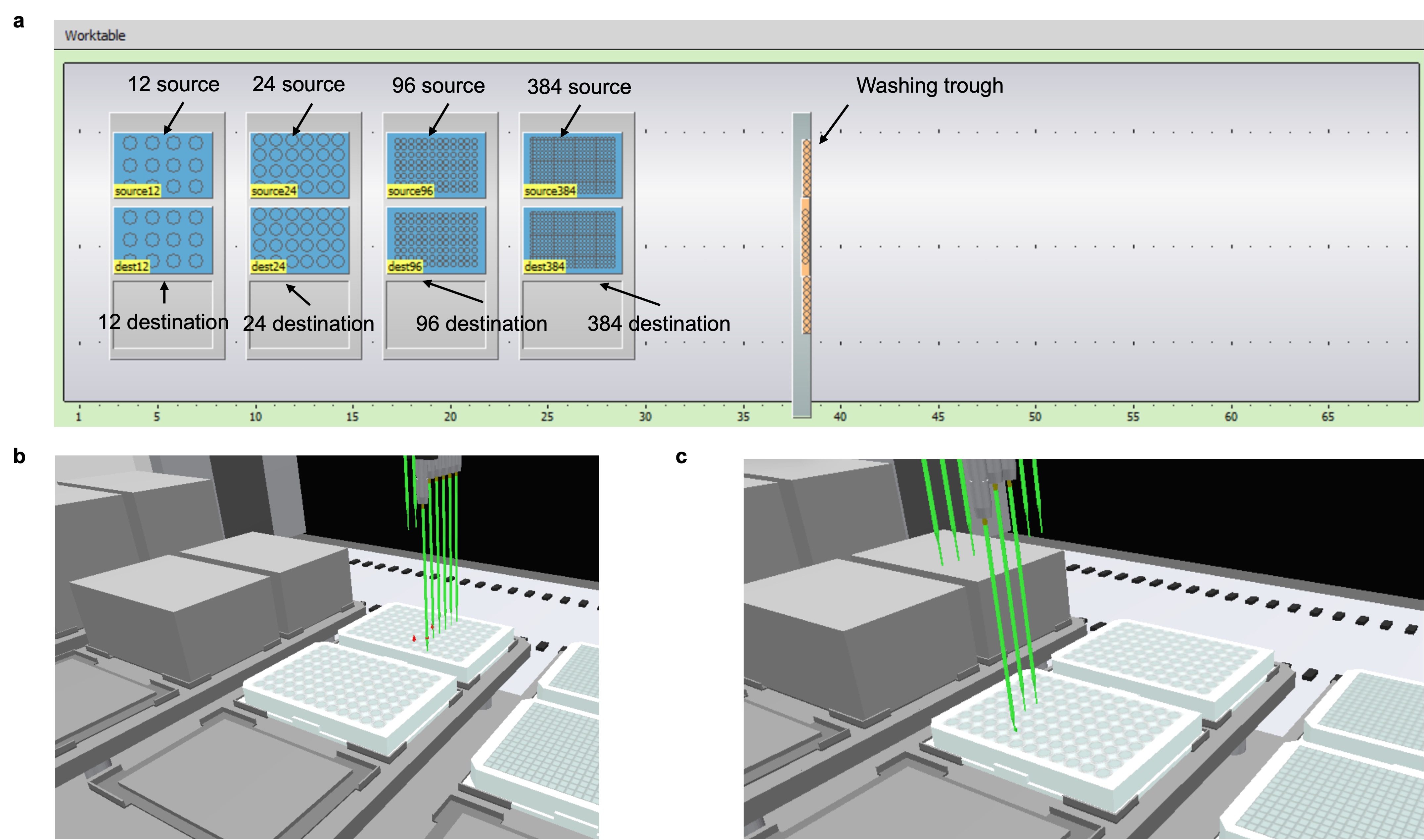}
\caption{(a)Worktable layout for simulation. (b) Snapshot of aspirating from 96-well plate. (c) Snapshot of depensing to 96-well plate.}
\label{demo_simulation}
\end{figure}

\begin{figure}[H]
\centering
\includegraphics[width=0.6\textwidth]{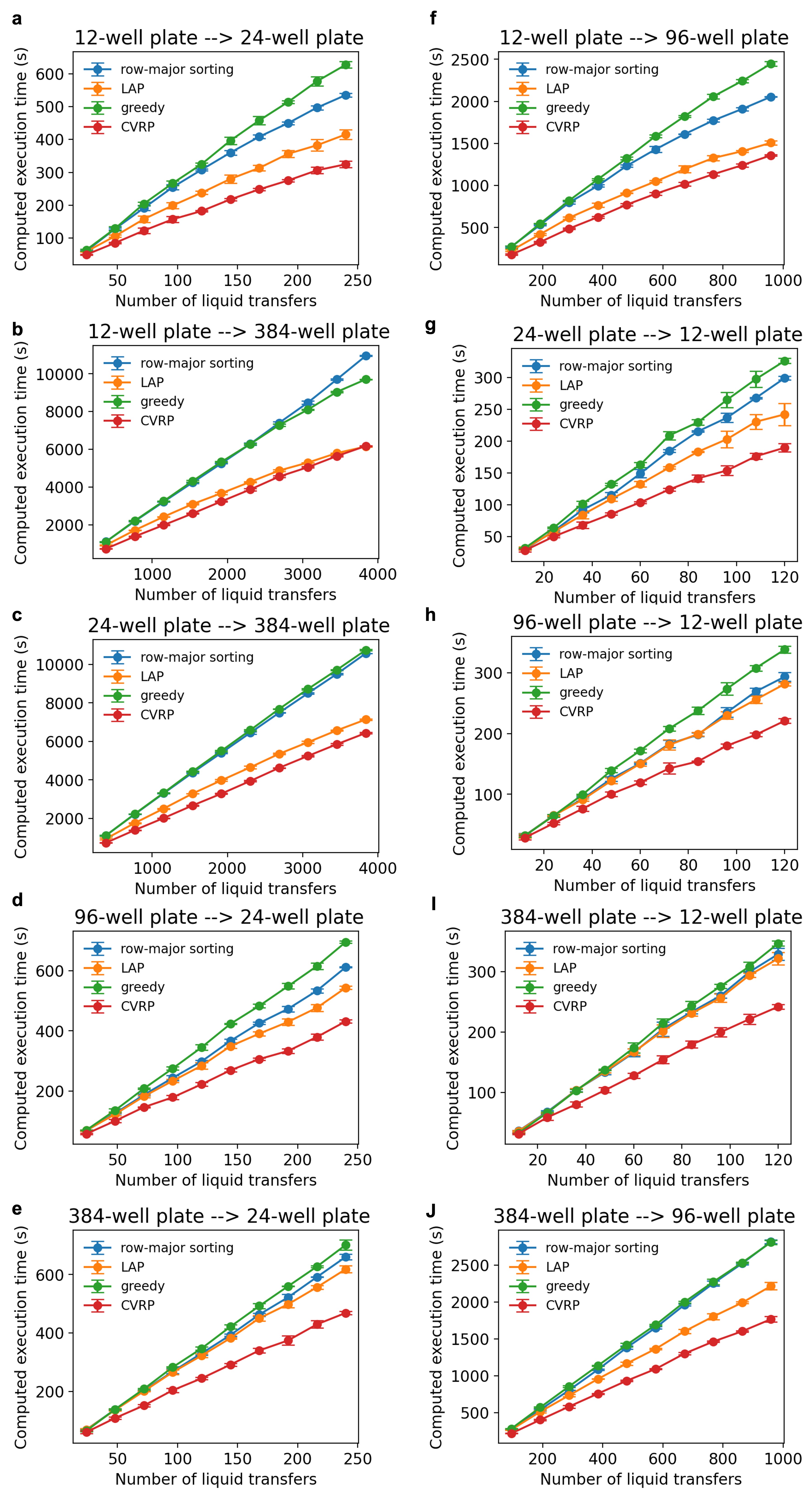}
\caption{Execution time across different labware formats for random tasks.(a) from 12-well plate to 24-well plate, (b) from 12-well plate to 384-well plate, (c) from 24-well plate to 384-well plate, (d) from 96-well plate to 24-well plate, (e) from 384-well plate to 24-well plate, (f) from 12-well plate to 96-well plate, (g) from 24-well plate to 12-well plate, (h) from 96-well plate to 12-well plate, (i) from 384-well plate to 12-well plate, and (j) from 384-well plate to 96-well plate. For a given number of liquid transfers and labwares, 3 random tasks were generated and scheduled using different methods. The solving time for CVRP is 20 seconds. For each labware combination, we evaluated 10 different task sizes, with the number of liquid transfers set to multiples (from 1× to 10×) of the number of wells in the destination plate.
}
\label{remaining_results}
\end{figure}

\begin{figure}[H]
\centering
\includegraphics[width=0.85\textwidth]{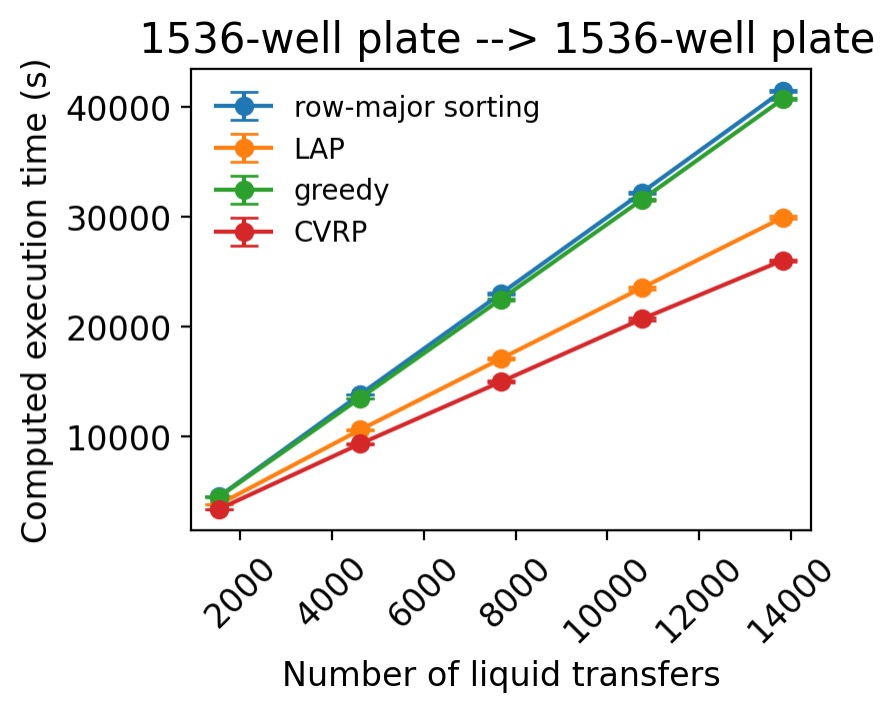}
\caption{Performance of different methods on 1536-well plates. For a given number of liquid transfers, 3 random tasks were generated and scheduled using different methods. The solving time for CVRP is 60 seconds for tasks with less than 7000 transfers and 120 seconds for tasks with more than 7000 transfers. We evaluated 5 different task sizes (1536, 4608, 7680, 10752 and 13824).}
\label{1536_results}
\end{figure}

\begin{figure}[H]
\centering
\includegraphics[width=0.85\textwidth]{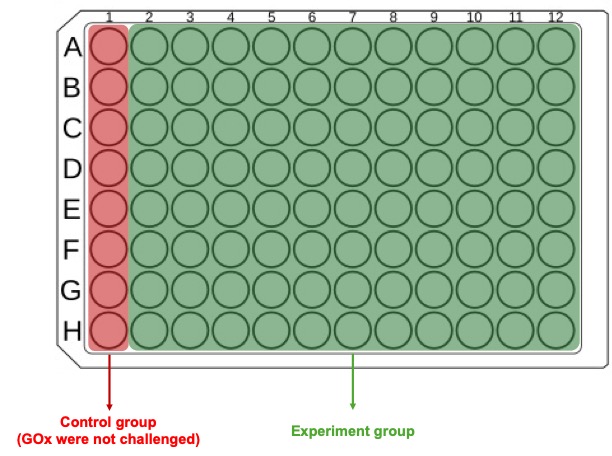}
\caption{The layout of the 96-well plate during the autonomous optimization of polymer blends for enzyme stability.}
\label{gox_layout}
\end{figure}

\begin{figure}[H]
\centering
\includegraphics[width=0.85\textwidth]{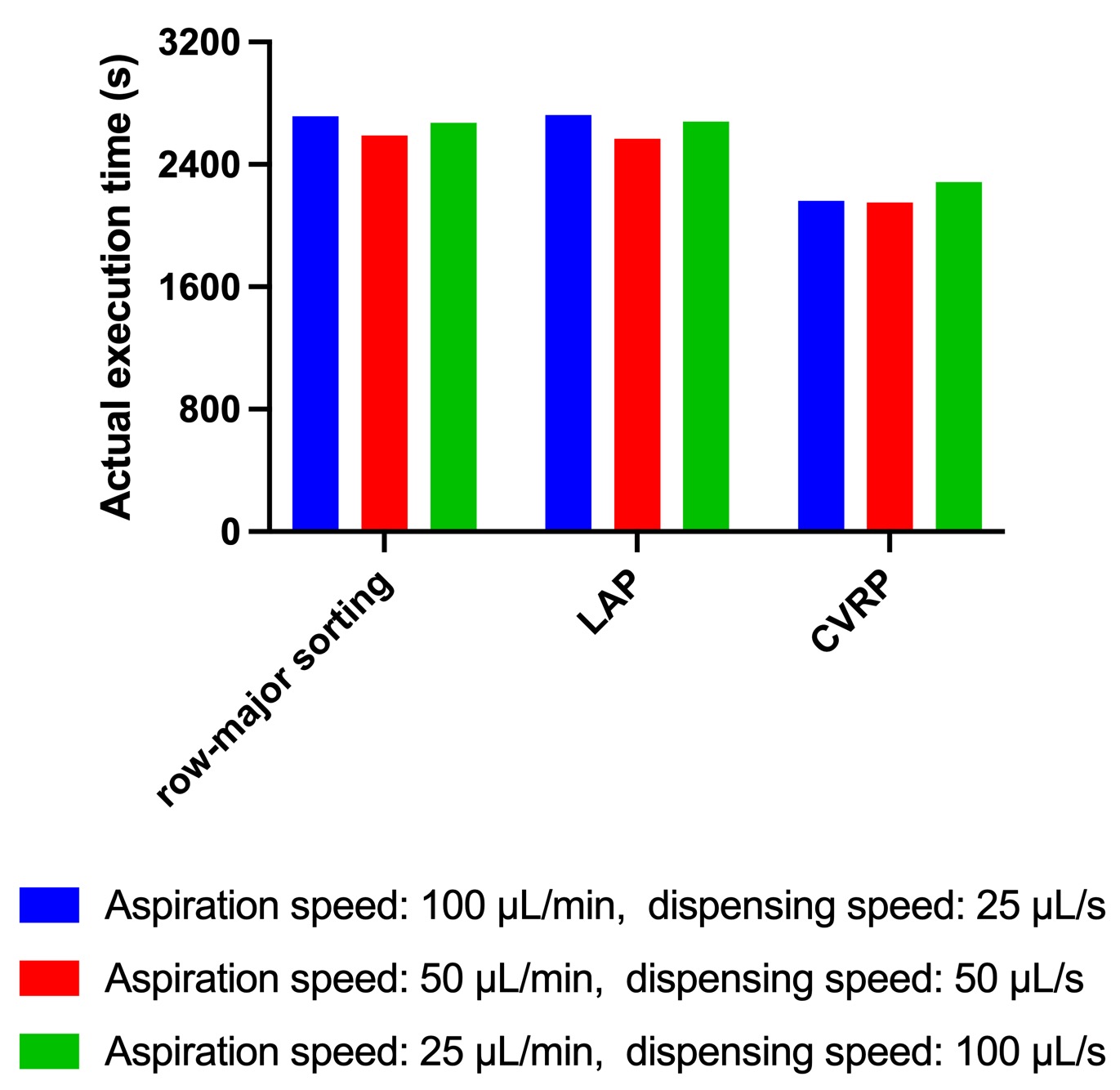}
\caption{The actual execution time on JANUS liquid handling workstation of iteration 3 in Figure 7 scheduled with different methods. The optimization effect was tested according to the parameters specified in Table \ref{janus_params}, with the specified aspiration and dispensing speeds. The CVRP-based method demonstrated superior performance compared to all other methods in various speed combinations. }
\label{janus_results}
\end{figure}